%
\documentstyle{amsppt}
\nopagenumbers
\font\eightuuu=cmr8
\font\twelveuuu=cmr12
\def\compos{\,\raise 1pt\hbox{$\sssize\circ$} \,}
\def\const{\operatorname{const}}
\def\Rea{\operatorname{Re}}
\def\Img{\operatorname{Im}}
\def\Jac{\operatorname{Jac}}
\def\Prym{\operatorname{Prym}}
\def\Res{\operatornamewithlimits{Res}}
\pagewidth{360pt}
\pageheight{606pt}
\loadbold
\TagsOnRight
\document
\leftline{\vbox{\eightuuu\hsize=3.5cm
\centerline{THEOR\. AND MATH\.}
\centerline{PHYSICS (RUSSIAN)}
\centerline{87(1991) No\. 1, 48--56}}}
\vskip 0.7cm
\leftline{{\twelveuuu\copyright} 1991\hskip 3.5cm {\bf R\.~A\.~Sharipov}}
\vskip 0.5cm
\centerline{\bf MINIMAL TORI IN FIVE-DIMENSIONAL SPHERE IN $\Bbb C^{\,3}$.}
\bigskip
\parshape 1 30pt 300pt {\bf Abstract.} Special class of surfaces in
five-dimensional sphere in $\Bbb C^{\,3}$ is considered. Immersion
equations for minimal tori of that class are shown to be reducible
to the equation $u_{z\bar z}=e^u-e^{-2u}$ which is integrable by
means of inverse scattering method. Finite-gap minimal tori are
constructed.\par
\head
1. Introduction.
\endhead
    Minimal surfaces in multidimensional spaces naturally arise
as classical trajectories of relativistic strings with Lagrangians
suggested by Nambu (see \cite{1}) and by Polyakov (see \cite{2}).
From geometric point of view these are surfaces of zero mean curvature.
Their immersion into environment space often is given by the equations
integrable by means of inverse scattering method (see \cite{3--5}).
Minimal surfaces in $\Bbb R^3$ and $\Bbb R^4$ are described by
Liouville equation $u_{z\bar z}=e^u$. This equation is nonlinear, but
can be linearized by Backlund transformation (see \cite{6, 7}).
In higher dimensional spaces and in non-flat spaces immersion equations
for minimal surfaces are not linearizable. But they possess Lax pairs
(see \cite{5}), therefore one can rather effectively study their
solutions.\par
    In this paper minimal tori in five-dimensional sphere $S^5\subset
\Bbb C^{\,3}$ are considered whose immersion is described by
Bullough-Dodd-Jiber-Shabat equation\footnote{After this paper had been
published S\.~ P\.~Tsarev discovered that the equation \thetag{1.1} was
first introduced by Tzitzeica in \cite{12}. Now it is called Tzitzeica
equation. See also \cite{13} for more details.}
$$
\hskip -2em
u_{z\bar z}=e^u-e^{-2u}.
\tag1.1
$$
On a base of construction of finite-gap solutions for this equation
(see \cite{8}) finite-gap minimal tori which are complexly normal in
$S^5$ are constructed. Situation here is similar to that considered
in \cite{9} and \cite{10}. There substantial advances in describing
tori of constant mean curvature in $\Bbb R^3$, in $S^3$, and in
$H^3$ on a base of finite gap solutions of Sine-Gordon equation
$u_{z\bar z}=\sin u$ were achieved. In the framework of affine geometry
the equation \thetag{1.1} were considered in \cite{11}.
\head
2. Complexly normal surfaces in Hermitian sphere in $\Bbb C^{\,3}$\\
and their scalar invariants.
\endhead
   Let's consider complex space $\Bbb C^{\,3}$ with standard Hermitian
scalar product
$$
\hskip -2em
\left<\bold A\,|\,\bold B\right>=\sum^3_{i=1}\bar A_i\,B_i
\tag2.1
$$
and with associated Euclidean scalar product
$$
\hskip -2em
\left(\bold A\,|\,\bold B\right)=\Rea\bigl(\left<\bold A\,|
\,\bold B\right>\bigr).
\tag2.2
$$
Let $\bold r(x^1,x^2)$ be vector-valued function with values in
$\Bbb C^{\,3}$ defining immersion of real two-dimensional surface
$T$ into the sphere $S_R$ of radius $R$ in $\Bbb C^{\,3}$. Denote
by $\bold E_1$ and $\bold E_2$ basic tangent vectors of this
two-dimensional surface:
$$
\xalignat 2
&\bold E_1=\partial_1\bold r=\frac{\partial\bold r}{\partial x^1},
&&\bold E_2=\partial_2\bold r=\frac{\partial\bold r}{\partial x^2}.
\endxalignat
$$
Scalar product \thetag{2.1} induces Hermitian metric on $T$:
$$
\hskip -2em
h_{ij}=\left<\bold E_i\,|\,\bold E_j\right>=g_{ij}+i\,\omega_{ij}.
\tag2.3
$$
Its real part is a Riemannian metric induced by scalar product
\thetag{2.2}, while imaginary part of Hermitian metric \thetag{2.3}
is skew-symmetric tensor
$$
\hskip -2em
\omega_{ij}=\left(i\cdot\bold E_i\,|\,\bold E_j\right)
\tag2.4
$$
defining closed $1$-form $\boldsymbol\omega$ on $T$. Operator-valued
tensor field
$$
\Omega^i_j=\sum^2_{k=1}g^{ik}\,\omega_{kj}
$$
has zero trace, while $\det\Omega$ is a scalar invariant of metric
\thetag{2.3} on $T$.
\definition{Definition} Immersion $\bold r\!: T\to M\subset\Bbb
C^{\,3}$ of real two-dimensional surface $T$ into real submanifold
$M$ of codimension $1$ in $\Bbb C^{\,3}$ is called {\it complexly
normal} immersion if at each point of $T$ Euclidean unit normal
vector $\bold N$ of $M$ is orthogonal to tangent plane to $T$
in Hermitian metric, i\.\,e\. $\left<\bold E_i\,|\,\bold N\right>=0$.
\enddefinition
    For complexly normal surface $T$ we define vectors $\bold F_1$
and $\bold F_2$ which are orthogonal to $\bold N$ in Hermitian metric
and orthogonal to vectors $\bold E_1$ and $\bold E_2$ in Euclidean
metric. We define them by the relationship
$$
\hskip -2em
\bold F_i=i\cdot\bold E_i+\sum^2_{s=1}\Omega^s_i\cdot E_s.
\tag2.5
$$
Vectors $\bold F_1$ and $\bold F_2$ define one more tensor field
related to $g_{ij}$ and $\omega_{ij}$:
$$
f_{ij}=g_{ij}+\sum^2_{k=1}\sum^2_{s=1}\omega_{ik}\,g^{ks}
\,\omega_{sj}.
$$
For associated tensor field $F^i_j$ we derive
$$
F^i_j=\sum^2_{k=1}g^{ik}\,f_{kj}=\delta^i_j+\sum^2_{k=1}
\sum^2_{r=1}\sum^2_{s=1}g^{ik}\,\omega_{kr}\,g^{rs}\,
\omega_{sj}=\delta^i_j+\sum^2_{r=1}\Omega^i_r\,\Omega^r_j.
$$
Scalar invariants of this field can be expressed through invariant
$\det\Omega$.\par
    Vectors $\bold E_1$ and $\bold E_2$ form a moving frame in tangent
space to surface $T$, while vectors $\bold F_1$, $\bold F_2$, $\bold N$,
and $i\cdot\bold N$ form complementary moving frame in normal space.
Dynamics of first frame is given by the equations
$$
\hskip -2em
\partial_i\bold E_j=\sum^2_{k=1}\Gamma^k_{ij}\cdot\bold E_k
+\sum^2_{k=1}T^k_{ij}\cdot\bold F_k+\left(b_{ij}+i\,d_{ij}
\right)\cdot\bold N,
\tag2.6
$$
where $\Gamma^k_{ij}$ are components of metric connection of $T$ given by
well-known formula
$$
\hskip -2em
\Gamma^k_{ij}=\sum^2_{s=1}\frac{g^{ks}}{2}\left(\partial_ig_{sj}+
\partial_jg_{is}-\partial_sg_{ij}\right).
\tag2.7
$$
Dynamics of vectors $\bold F_1$ and $\bold F_2$ is completely determined
by dynamics of $\bold E_1$ and $\bold E_2$ due to the relationship
\thetag{2.5}. Dynamics of unit normal vector $\bold N$ of submanifold
$M$ along $T$ is given by the following equation:
$$
\hskip -2em
\partial_i\bold N=\sum^2_{k=1}L^k_i\cdot\bold E_k+\sum^2_{k=1}M^k_i
\cdot\bold F_k+i\,S_i\cdot\bold N.
\tag2.8
$$
Tensor fields $b_{ij}$, $d_{ij}$, $L^k_i$, and $M^k_i$, which appear as
coefficients in the equations \thetag{2.6} and \thetag{2.8}, are related
to each other by a number of relationships that follow from our specific
choice of frames. Due to orthogonality of $\bold N$ and $\bold E_j$ we
get
$$
\hskip -2em
b_{ij}=-\sum^2_{k=1}L^k_i\,g_{kj},
\tag2.9
$$
while orthogonality of vectors $i\cdot\bold N$ and $\bold E_j$ yields
$$
\hskip -2em
d_{ij}=\sum^2_{k=1}M^k_i\,f_{kj}-\sum^2_{k=1}L^k_i\,\omega_{kj}.
\tag2.10
$$
Differentiating \thetag{2.4} and using the relationship \thetag{2.6},
we derive the equality
$$
\hskip -2em
\nabla_s\omega_{ij}=\sum^2_{k=1}T^k_{sj}\,f_{ki}
-\sum^2_{k=1}T^k_{si}\,f_{kj},
\tag2.11
$$
which completely determines skew-symmetric in indices $i$ and $j$
part of tensor
$$
t_{isj}=\sum^2_{k=1}T^k_{is}\,f_{kj}.
$$\par
    In the case, when manifold $M$ is a sphere $S_R$ of radius $R$, the
above relationships \thetag{2.8}--\thetag{2.11} simplify substantially.
In this case radius-vector $\bold r(x^1,x^2)$ is collinear to normal
vector: $\bold r=R\cdot\bold N$. Therefore
$$
\bold E_i=\partial_i\bold r=R\cdot\partial_i\bold N.
$$
Comparing this equality with \thetag{2.8}, we get
$$
\xalignat 3
&L^k_i=\frac{1}{R}\,\delta^k_i,
&&M^k_i=0,
&&S_i=0.
\endxalignat
$$
Further from \thetag{2.9} and \thetag{2.10} for matrices of second
fundamental forms we derive
$$
\xalignat 2
&b_{ij}=-\frac{1}{R}\,g_{ij},
&&d_{ij}=-\frac{1}{R}\,\omega_{ij}.
\endxalignat
$$
Matrix $d_{ij}$ is symmetric, while matrix $\omega_{ij}$ is
skew-symmetric. Therefore both matrices are zero: $d_{ij}=\omega_{ij}
=0$. Thus, for $M=S_R$ and for complexly normal embedding of surface
$T$ vectors $\bold F_1$ and $\bold F_2$ coincide with vectors $i\cdot
\bold E_1$ and $i\cdot\bold E_2$ respectively, while relationships
\thetag{2.6} and \thetag{2.8} are written as
$$
\xalignat 2
&\hskip -2em
\nabla_i\bold E_j=\sum^2_{k=1}T^k_{ij}\cdot\bold F_k-\frac{1}{R}
\,g_{ij}\cdot\bold N,
&&\partial_i\bold N=\frac{1}{R}\,\bold E_i.
\tag2.12
\endxalignat
$$
Due to the equality \thetag{2.11} tensor
$$
T_{kij}=\sum^2_{s=1}T^s_{ij}\,g_{sk}
$$
is symmetric with respect to all its indices. Gauss equation,
Peterson-Coddazi equation and Ricci equation are obtained as compatibility
conditions for the equations \thetag{2.12}. Here is Gauss equation
$$
\hskip -2em
R^s_{kij}=\sum^2_{r=1}T^r_{jk}\,T^s_{ir}-\sum^2_{r=1}T^r_{ik}\,T^s_{jr}
+\frac{g_{jk}\,\delta^s_i-g_{ik}\,\delta^s_j}{R^2},
\tag2.13
$$
where $R^s_{kij}$ is Riemann curvature tensor determined by metric
connection \thetag{2.7} according to standard formula
$$
\hskip -2em
R^s_{kij}=\partial_i\Gamma^s_{kj}-\partial_j\Gamma^s_{ki}-
\sum^2_{r=1}\Gamma^r_{ki}\,\Gamma^s_{rj}+\sum^2_{r=1}\Gamma^r_{kj}
\,\Gamma^s_{ri}.
\tag2.14
$$
Peterson-Coddazi and Ricci equations in this case are united into one
equation
$$
\hskip -2em
\nabla_iT_{jsk}-\nabla_jT_{isk}=0.
\tag2.15
$$
Symmetric tensor $T_{ijk}$ has two scalar invariants of second order
$$
\xalignat 2
&\hskip -2em
H^2=\sum^2_{i=1}\sum^2_{j=1}\sum^2_{s=1}T^{is}_i\,T^j_{js},
&&k=\sum^2_{i=1}\sum^2_{j=1}\sum^2_{s=1}T^{ijs}\,T_{ijs}
\tag2.16
\endxalignat
$$
and an invariant of fourth order determined by the relationship
$$
\hskip -2em
q=\sum^2_{i=1}\sum^2_{j=1}\sum^2_{k=1}\sum^2_{r=1}\sum^2_{p=1}
\sum^2_{s=1}T^i_{jk}\,T^{jk}_s\,T^s_{rp}\,T^{rp}_i.
\tag2.17
$$
Due to specific features of two-dimensional case ($\dim T=2$)
invariants \thetag{2.16} and \thetag{2.17} form maximal set
of functionally independent invariants of symmetric tensor
$T_{ijk}$. Moreover, in two-dimensional case due to symmetry
of Riemann curvature tensor $R^s_{kij}$ tensorial equation
\thetag{2.13} is equivalent to one scalar equation that binds
Gaussian curvature $K$ of the surface $T$ with curvatures $H$
and $k$ of tensor $T_{ijk}$:
$$
\hskip -2em
2\,K=\sum^2_{j=1}g^{kj}\,R^s_{ksj}=H^2-k+2\,R^{-2}.
\tag2.18
$$
Invariant $H$ in \thetag{2.16}
coincides with the length of averaged normal vector of $T$:
$$
\pagebreak
\hskip -2em
H\cdot\bold n=\sum^2_{i=1}T^{ik}_i\cdot\bold F_k.
\tag2.19
$$
Here $H$ is mean curvature of the surface $T$ embedded into $S$,
while unit vector $\bold n$ tangent to sphere $S_R$ in \thetag{2.19}
is a unit vector of averaged normal of $T$.
\head
3. Complexly normal tori of zero mean curvature.
\endhead
    The condition of vanishing of mean curvature is very restricting
condition for the class of surfaces in question. Indeed, from the
condition $H=0$ due to the equality \thetag{2.19} we have
$$
\sum^2_{i=1}T^{ik}_i=0.
$$
If we take into account symmetry of tensor $T_{ijk}$, the above equality
means that in this tensor we have only two independent components. In
order to use this circumstance let's choose isothermal coordinates on the
surface $T$, i\.\,e\. coordinates $x=x^1=\Rea z$ and $y=x^2=\Img z$ for
which metric $g_{ij}$ is conformally Euclidean: $\bold g=2\,R^2\,e^u\,
dz\,d\bar z$. In this case components of tensor $T^k_{ij}$ are expressed
through two independent quantities $A$ and $B$:
$$
\xalignat 2
&\hskip -2em
T^1_{11}=A,&&T^2_{12}=T^2_{21}=T^1_{22}=-A,\\
\vspace{-1.5ex}
&&&\tag3.1\\
\vspace{-1.5ex}
&\hskip -2em
T^2_{22}=B,&&T^1_{12}=T^1_{21}=T^2_{11}=-B.
\endxalignat
$$
Let's calculate components of metric connection $\Gamma^k_{ij}$ by using
formula \thetag{2.7}:
$$
\xalignat 3
&\hskip -2em
\Gamma^1_{11}=\frac{u_x}{2},&&\Gamma^2_{11}=-\frac{u_y}{2},
&&\Gamma^1_{12}=\Gamma^1_{21}=\frac{u_y}{2},\\
\vspace{-1.3ex}
&&&\tag3.2\\
\vspace{-1.3ex}
&\hskip -2em
\Gamma^2_{22}=\frac{u_y}{2},&&\Gamma^1_{22}=-\frac{u_x}{2},
&&\Gamma^2_{12}=\Gamma^2_{21}=\frac{u_x}{2}.
\endxalignat
$$
Then let's substitute \thetag{3.1} into Peterson-Coddazi-Ricci
equation \thetag{2.15}. Upon completing calculations and
taking into account \thetag{3.2} we find
$$
\xalignat 2
&\hskip -2em
\partial_x\left(e^u\,A\right)=\partial_y\left(e^u\,B\right),
&&\partial_y\left(e^u\,A\right)=-\partial_x\left(e^u\,B\right).
\tag3.3
\endxalignat
$$
It's easy to see that relationships \thetag{3.3} do coincide with
Cauchy-Riemann equations for holomorphic function $G(z)=e^u\,A+
i\,e^u\,B$.\par
    If $G(z)$ is identically zero, then we have trivial case. In this
case due to \thetag{2.12} the subspace defined as a span of vectors
$\bold E_1$, $\bold E_2$, and $\bold N$ contains all derivatives of
these vectors. Hence this subspace do not change when we vary $x$ and
$y$. This means that vectors $\bold E_1$, $\bold E_2$, and $\bold N$
belong to some fixed three-dimensional real subspace of $\Bbb C^{\,3}$,
while $T$ is a central section of the sphere $S_R$ within this
three-dimensional subspace, i\.\,e\. $T$ is two-dimensional sphere
of radius $R$ or its part.\par
   Let's consider the case $G(z)\not\equiv 0$. In this case we shall
assume that $T$ is compact closed surface of toric topology. Taking
$z$ for uniformizing parameter inherited from universal cover $\Bbb C
\to T$, we get $G(z)=\const\neq 0$ since $G(z)$ then is holomorphic
function on complex torus $T$. At the expense of simultaneous change
of scale along axes $x$ and $y$ by the same factor we can satisfy
additional \pagebreak condition $|G(z)|=1$. Then we write $G(z)=\cos
\vartheta+i\,\sin\vartheta$. Now \thetag{3.1} is written as
$$
\xalignat 2
&\hskip -2em
T^1_{11}=e^{-u}\,\cos\vartheta,
&&T^2_{22}=e^{-u}\,\sin\vartheta,\\
&\hskip -2em
T^2_{11}=-e^{-u}\,\sin\vartheta,
&&T^1_{22}=-e^{-u}\,\cos\vartheta,
\tag3.4\\
&\hskip -2em
T^1_{12}=T^1_{21}=-e^{-u}\,\sin\vartheta,
&&T^2_{12}=T^2_{21}=e^{-u}\,\cos\vartheta.
\endxalignat
$$
Using tensor $T^k_{ij}$ of the form \thetag{3.4}, we can calculate
invariants $k$ and $q$ determined by formulas \thetag{2.16} and
\thetag{2.17}:
$$
\xalignat 2
&\hskip -2em
k=\frac{2\,e^{-3u}}{R^2},
&&q=\frac{2\,e^{-6u}}{R^4}.
\tag3.5
\endxalignat
$$
If we calculate curvature tensor $R^s_{kij}$ for connection components
\thetag{3.2} using formula \thetag{2.14} and if we substitute it into
Gauss equation \thetag{2.18}, we obtain the equation for the function
$u(x,y)$:
$$
u_{xx}+u_{yy}=4\,e^{-2u}-4\,e^u,
$$
This equation does coincide with \thetag{1.1}. Its solution corresponding
to the embedding of two-dimensional torus into $S_R\subset\Bbb C^{\,3}$ is
double-periodic with some grid of periods in the plane of variables $x$
and $y$. Below we consider class of finite-gap minimal surfaces in sphere
$S_R$ including compact minimal tori which are complexly normal in this
sphere.
\head
Finite gap solutions of the equation $u_{z\bar z}=e^u-e^{-2u}$ and
associated orthonormal frame in $\Bbb C^{\,3}$.
\endhead
    Let's consider Riemannian surface $\Gamma$ of even genus $g$ with
two distinguished points $P_0$ and $P_\infty$ such that
\roster
\rosteritemwd=0pt
\item"1)" there is a meromorphic function $\lambda(P)$ on $\Gamma$ with
divisor of zeros and poles $3\,P_0-3\,P_\infty$;
\item"2)" there is a holomorphic involution $\sigma$ such that
$\lambda(\sigma P)=-\lambda(P)$;
\item"2)" there is an antiholomorphic involution $\tau$ of separating
type such that
$$
\hskip -2em
\lambda(\tau P)\overline{\lambda(P)}=1.
\tag4.1
$$
\endroster
Points of $\Gamma$ which are stable under the action of $\tau$
(i\.\,e\. $\tau P=P$) form a closed curve or a set of several
closed curves. These curves are called invariant cycles of
$\tau$. For antiholomorphic involution of separating type $\tau$
its invariant cycles break $\Gamma$ into two domains: $\Gamma_0$
containing point $P_0$ and $\Gamma_\infty$ containing point
$P_\infty$. Due to \thetag{4.1} all invariant cycles of $\tau$
are projected onto unit circle on complex $\lambda$-plane. The
number of these cycles is less or equal to three. This number
determines the number of real tori in Jacobian $\Jac(\Gamma)$.
Each such torus is composed by classes of divisors $D$ of degree
$g$ such that
$$
D+\tau D-P_0-P_\infty=C,
$$
where $C$ is divisor of canonic class on $\Gamma$. Due to
\thetag{4.2} each real divisor $D$ (e\.\,e\. divisor from
real torus in $\Jac(\Gamma)$) determines some Abelian differential
of the third kind $\omega(P)$ with zeros at the points of divisor
$D+\tau D$ and with residues
$$
\xalignat 2
&\hskip -2em
\Res_{P=P_0}\omega(P)=+i,
&&\Res_{P=P_\infty}\omega(P)=-i,
\endxalignat
$$
at the points $P_0$ and $P_\infty$, where $\omega(P)$ has simple poles.
Under the action of $\tau$ differential $\omega(P)$ is transformed as
follows:
$$
\hskip -2em
\omega(\tau P)=\overline{\omega(P)}.
\tag4.3
$$
Therefore it is real valued on invariant cycles of $\tau$. Real torus
$T_0$ in $\Jac(\Gamma)$ is distinguished among other real tori by the
following property: for divisor $D$ from this torus differential
$\omega(P)$ is positive  on all invariant cycles of $\tau$ with respect
to natural orientation of boundary $\partial\Gamma_\infty$.\par
    Upon fixing torus $T_0$ in $\Jac(\Gamma)$ let's consider its subset
consisting of divisors invariant under the action of composite map
$\tau\compos\sigma$:
$$
\hskip -2em
\tau D=\sigma D.
\tag4.4
$$
This subset is not empty. It is real torus $T_0$ in Prym variety
$\Prym(\Gamma)$. For divisors of this torus we can complete
\thetag{4.3} by another relationship
$$
\omega(\sigma P)=\omega(P),
$$
which follows from \thetag{4.4} and from invariance of points $P_0$
and $P_\infty$ under the action of involution $\sigma$.\par
    Let's fix local parameters $k^{-1}(P)$ and $q^{-1}(P)$ in the
neighborhood of distinguished points $P_0$ and $P_\infty$ by the
conditions
$$
\xalignat 2
&\hskip -2em
k^3(P)=\lambda(P),
&&\overline{k(\tau P)}=q(P).
\tag4.5
\endxalignat
$$
Now, having fixed some positive divisor $D\in T_0\subset\Prym(\Gamma)$
of degree $g$, we construct vectorial Baker-Achiezer function
$\psi(z,P)$ with values in $\Bbb C^{\,3}$ such that
$$
\hskip -2em
\aligned
&\psi_1(P)=e^{i\,k(P)\,z}\left(k^{-1}(P)+\ldots\right),\\
&\psi_2(P)=e^{i\,k(P)\,z}\left(k^{-2}(P)+\ldots\right),\\
&\psi_1(P)=e^{i\,k(P)\,z}\left(k^{-3}(P)+\ldots\right)
\endaligned
\tag4.6
$$
in the neighborhood of distinguished point $P_\infty$ and such that
$$
\hskip -2em
\aligned
&\psi_1(P)=e^{i\,q(P)\,\bar z}\left(q^{1}(P)+\ldots\right)e^{-u},\\
&\psi_2(P)=e^{i\,q(P)\,\bar z}\left(q^{2}(P)+\ldots\right)e^u,\\
&\psi_1(P)=e^{i\,q(P)\,\bar z}\left(q^{3}(P)+\ldots\right)
\endaligned
\tag4.7
$$
in the neighborhood of another distinguished point $P_0$. Functions
$\psi_1$, $\psi_2$, and $\psi_3$ are uniquely determined by divisor
$D$ and by conditions \thetag{4.5} and \thetag{4.7} (see \cite{8}).
They satisfy the following differential equations:
$$
\xalignat 2
&\hskip -2em
\partial_z\psi_1=-u_z\,\psi_1+i\,\lambda\,\psi_3,
&&\partial_{\bar z}\psi_1=i\,e^{-2u}\,\psi_2,\\
&\hskip -2em
\partial_z\psi_2=u_z\,\psi_2+i\,\psi_1,
&&\partial_{\bar z}\psi_2=i\,e^{u}\,\psi_3,
\tag4.8\\
&\hskip -2em
\partial_z\psi_3=i\,\psi_2,
&&\partial_{\bar z}\psi_2=i\,\lambda^{-1}\,e^u\,\psi_1.\\
\endxalignat
$$
Compatibility condition of these equations is equivalent to the
equation \thetag{1.1} for the function $u=u(z,\bar z)$ in
\thetag{4.8}. Condition $D\in T_0\subset\Prym(\Gamma)$ provides
that $u$ is real-valued smooth function. There is an explicit
formula for $u$ in terms of Prym theta-functions (see \cite{8}).
Respective to $\psi_1$, $\psi_2$, and $\psi_3$ the same condition
expressed by \thetag{4.4} and \thetag{4.5} yields
$$
\hskip -2em
\aligned
&\psi_1(\sigma P)=-\lambda^{-1}(P)\,e^{-u}\,\overline{\psi_2(\tau P)},\\
&\psi_2(\sigma P)=\lambda^{-1}(P)\,e^{u}\,\overline{\psi_1(\tau P)},\\
&\psi_3(\sigma P)=-\lambda^{-2}(P)\,\overline{\psi_2(\tau P)}.
\endaligned
\tag4.9
$$\par
     Remarkable feature of spectral problems associated with integrable
nonlinear equations is the presence of bilinear forms (pairings or
generalized Wronskians), which are in concordance with Lax operators
and which in finite-gap case possess some ``resonant'' properties. The
latter property can be used in constructing soliton-like solutions
on finite-gap background for these equations and in Cauchy kernels on
Riemann surfaces. For spectral problem \thetag{4.8} such pairing is given
by formula
$$
\hskip -2em
\gathered
\Omega(P,Q)=\left\{\psi(P)\,|\,\psi(\sigma Q)\right\}=
\psi_1(P)\,\psi_2(\sigma Q)\,\lambda(P)\,-\\
-\,\psi_2(P)\,\psi_1(\sigma Q)\,\lambda(P)-
\psi_3(P)\,\psi_3(\sigma Q)\,\lambda^2(P).
\endgathered
\tag4.10
$$
Let's differentiate \thetag{4.10} with respect to $z$ and $\bar z$.
Taking into account differential equations \thetag{4.8}, we obtain
the equalities
$$
\align
&\partial_z\Omega(P,Q)=i\left(\lambda(Q)-\lambda(P)\right)
\lambda(P)\,\psi_2(P)\,\psi_3(\sigma Q),\\
&\partial_{\bar z}\Omega(P,Q)=i\,e^u\,\left(\lambda(P)
\,\lambda^{-1}(Q)-1\right)\lambda(P)\,\psi_3(P)\,\psi_1(\sigma Q)
\endalign
$$
Looking at these equalities, we see that for $Q=P$ function \thetag{4.10}
does not depend on $z$ and $\bar z$. Moreover, function $W(P)=
\Omega(P,P)$ is meromorphic on $\Gamma$, it can be calculated explicitly
in the following form:
$$
\hskip -2em
W(P)=\frac{i\,d\lambda(P)}{\lambda(P)\,\omega(P)}.
\tag4.11
$$
Note that $\lambda\!:\Gamma\to\Bbb C$ is three-sheeted covering.
Therefore each value $\lambda$ of the function $\lambda(P)$ has
multiplicity $3$, i\.\,e\. in general case there are three
distinct points $P_1$, $P_2$, and $P_3$ such that $\lambda(P_1)
=\lambda(P_2)=\lambda(P_3)=\lambda$. Resonant property of
$\Omega(P,Q)$ then is expressed by the following equality:
$$
\hskip -2em
\Omega(P_i,P_j)=\cases W(P_i)&\text{for \ }P_i=P_j,\\
\ \ 0&\text{for \ }P_i\neq P_j.\endcases
\tag4.12
$$
For each value of $\lambda$ such that $|\lambda|=1$ points
$P_1$, $P_2$, and $P_3$ are on invariant cycles of antiholomorphic
involution $\tau$. They are stable under the action of $\tau$.
Using them, we can compose the matrix $U=U(\lambda,z,
\bar z)$ of the following form:
$$
\hskip -2em
U=\Vmatrix
\dfrac{e^{u/2}\,\psi_1(P_1)}{\sqrt{W(P_1)}} &
\dfrac{e^{-u/2}\,\psi_2(P_1)}{\sqrt{W(P_1)}} &
\dfrac{\psi_3(P_1)}{\sqrt{W(P_1)}}\\
\vspace{1ex}
\dfrac{e^{u/2}\,\psi_1(P_2)}{\sqrt{W(P_2)}} &
\dfrac{e^{-u/2}\,\psi_2(P_2)}{\sqrt{W(P_2)}} &
\dfrac{\psi_3(P_2)}{\sqrt{W(P_2)}}\\
\vspace{1ex}
\dfrac{e^{u/2}\,\psi_1(P_3)}{\sqrt{W(P_3)}} &
\dfrac{e^{-u/2}\,\psi_2(P_3)}{\sqrt{W(P_3)}} &
\dfrac{\psi_3(P_3)}{\sqrt{W(P_1)}}
\endVmatrix
\tag4.13
$$
Resonant property \thetag{4.12}, invariance of $P_1$, $P_2$, and $P_3$
under the action of $\tau$, and the equalities \thetag{4.9} lead to
the following relationship:
$$
e^u\,\psi_1(P_i)\,\overline{\psi_1(P_j)}+
e^{-u}\,\psi_2(P_i)\,\overline{\psi_2(P_j)}+
\psi_3(P_i)\,\overline{\psi_3(P_j)}=W(P_i)\,\delta_{ij}.
$$
This means that matrix $U$ in \thetag{4.13} is unitary matrix. Moreover,
this equality means that function \thetag{4.11} is real and non-negative
on invariant cycles of $\tau$. Therefore square roots in \thetag{4.13}
are real numbers. Columns of unitary matrix $U$ form an orthonormal
frame in $\Bbb C^{\,3}$:
$$
\xalignat 3
&\bold L=U_1, &&\bold M=U_2, &&\bold N=U_3.
\tag4.14
\endxalignat
$$
This frame consists of three unit vectors perpendicular to each other
with respect to Hermitian metric \thetag{2.1}.
\head
5. Finite-gap embeddings of two-dimensional surfaces in $\Bbb C^{\,3}$.
\endhead
    Let's study the dynamics of orthonormal frame \thetag{4.14}. In
complex variable $z$ and $\bar z$ it is determined by the equations
\thetag{4.8}:
$$
\xalignat 2
&\partial_z\bold L=-\frac{u_z}{2}\cdot\bold L+i\,\lambda\,e^{u/2}
\cdot\bold N,
&&\partial_{\bar z}\bold L=\frac{u_{\bar z}}{2}\cdot\bold L
+i\,e^{-u}\,\cdot\bold M,\qquad\\
\vspace{1ex}
&\partial_z\bold M=\frac{u_z}{2}\cdot\bold M+i\,e^{u}
\cdot\bold L,
&&\partial_{\bar z}\bold M=-\frac{u_{\bar z}}{2}\cdot\bold M
+i\,e^{u/2}\,\cdot\bold N,\qquad
\tag5.1\\
\vspace{1ex}
&\partial_z\bold N=i\,e^{u/2}\cdot\bold M,
&&\partial_{\bar z}\bold N=i\,\lambda^{-1}\,e^{u/2}\,\cdot\bold L.
\endxalignat
$$
Passing to real variables $x=x^1=\Rea z$ and $y=x^2=\Img z$, from
\thetag{5.1} we derive the following equalities for the
dynamics of frame \thetag{4.14} with respect to $x$:
$$
\hskip -2em
\aligned
&\partial_x\bold L=\frac{u_y}{2}\cdot\bold L+i\,\lambda\,e^{u/2}
\cdot\bold N+i\,e^{-u}\cdot\bold M,\\
\vspace{1ex}
&\partial_x\bold M=-i\,\frac{u_y}{2}\cdot\bold M+i\,e^{-u}
\cdot\bold L+i\,e^{u/2}\cdot\bold N,\\
\vspace{1ex}
&\partial_x\bold N=i\,e^{u/2}\cdot\bold M+i\,\lambda^{-1}\,e^{u/2}\cdot
\bold L,
\endaligned
\tag5.2
$$
Similar  equations describe the dynamics of this frame with respect
to $y$:
$$
\hskip -2em
\aligned
&\partial_y\bold L=-\frac{u_x}{2}\cdot\bold L-\lambda\,e^{u/2}
\cdot\bold N+e^{-u}\cdot\bold M,\\
\vspace{1ex}
&\partial_y\bold M=i\,\frac{u_x}{2}\cdot\bold M-e^{u}
\cdot\bold L+e^{u/2}\cdot\bold N,\\
\vspace{1ex}
&\partial_y\bold N=-e^{u/2}\cdot\bold M+\lambda^{-1}\,e^{u/2}
\cdot\bold L.
\endaligned
\tag5.3
$$
Let's define the embedding of the surface $T$ into the sphere
$S_R\subset\Bbb C^{\,3}$ parametrically by means of function
$$
\hskip -2em
\bold r(x^1,x^2)=R\cdot\bold N(x^1,x^2)=R\cdot\bold N(x,y).
\tag5.4
$$
For tangent vectors $\bold E_1$ and $\bold E_2$ in the case of such
embedding we obtain
$$
\hskip -2em
\aligned
&\bold E_1=i\,R\,e^{u/2}\cdot\bold M+i\,R\,\lambda^{-1}\,e^{u/2}\cdot
\bold L,\\
\vspace{1ex}
&\bold E_2=-R\,e^{u/2}\cdot\bold M+R\,\lambda^{-1}\,e^{u/2}
\cdot\bold L.
\endaligned
\tag5.5
$$
Using Hermitian orthogonality of frame \thetag{4.14}, one can easily
find that vectors $\bold E_1$ and $\bold E_2$ are both orthogonal to
unit vector $\bold N$ with respect to Hermitian metric \thetag{2.1}.
Hence embedding \thetag{5.4} is complexly normal. Moreover, metric
$g_{ij}$ determined by \thetag{2.3} is diagonal and has conformally
Euclidean form $\bold g=2\,R^2\,e^u\left(dx^2+dy^2\right)$.\par
    Now let's apply the equations \thetag{5.2} and \thetag{5.3}.
Remember that $\lambda$ is chosen to be complex number of unit
modulus. Taking $\lambda=\cos\vartheta+i\,\sin\vartheta$, from
\thetag{5.5} we derive the dynamics of vectors $\bold E_1$ and
$\bold E_2$. Here it is:
$$
\hskip 2em
\aligned
&\nabla_1\bold E_1=e^{-u}\,\cos\vartheta\cdot\bold F_1-
e^{-u}\,\sin\vartheta\cdot\bold F_2-2\,R\,e^u\cdot\bold N,\\
\vspace{1ex}
&\nabla_2\bold E_1=\nabla_1\bold E_2=-e^{-u}\,\sin\vartheta\cdot
\bold F_1-e^{-u}\,\cos\vartheta\cdot\bold F_2,\\
\vspace{1ex}
&\nabla_2\bold E_2=-e^{-u}\,\cos\vartheta\cdot\bold F_1+
e^{-u}\,\sin\vartheta\cdot\bold F_2-2\,R\,e^u\cdot\bold N.
\endaligned
\tag5.6
$$
Components of metric connection for covariant derivatives in
\thetag{5.6} are given by \thetag{3.2}. Comparing \thetag{5.6}
with \thetag{2.12} we obtain components of tensor $T^k_{ij}$ for
the embedding \thetag{5.4}. They have exactly the same form as
given by \thetag{3.4}. Scalar invariants of this tensor are
given by formulas \thetag{3.5}. Compact finite-gap tori arise
in the case of double-periodic function $\psi_2$. Problem of
finding periodic solutions is standard in the theory of
finite-gap integration. They are obtained by imposing some
rather non-explicit restrictions to Riemann surface $\Gamma$,
which are written as rationality condition for some quotients
of Abelian integrals on it.

\enddocument
\end